\documentclass[11pt,amsmath,amssym,latexsymb]{article}
\usepackage{amssymb,latexsym}
\usepackage{longtable,tabularx}
\usepackage{makeidx}
\usepackage{graphicx}
\usepackage{psfrag}
\begin {document}

\def\newtext{}

\def\hat{\widehat}
\def\tilde{\widetilde}
\def \bfo {\begin {eqnarray*} }
\def \efo {\end {eqnarray*} }
\def \ba {\begin {eqnarray*} }
\def \ea {\end {eqnarray*} }
\def \beq {\begin {eqnarray}}
\def \eeq {\end {eqnarray}}
\newtheorem {lemma}{Lemma}\newtheorem {proposition}{Proposition}
\newtheorem {theorem}{Theorem}
\newtheorem {definition}{Definition}

\newtheorem {corollary}{Corollary}

\def\picture #1 by #2 (#3){
  \vsquare to #2{
    \hrule width #1 height 0pt depth 0pt
    \hfill
    \special{picture #3}}}

\def\scaledpicture #1 by #2 (#3 scaled #4){{
  \dimen0=#1 \dimen1=#2
  \divide\dimen0 by 1000 \multiply\dimen0 by #4
  \divide\dimen1 by 1000 \multiply\dimen1 by #4
  \picture \dimen0 by \dimen1 (#3 scaled #4)}}

\def \cM {{\cal M}}

\def \D {{\cal D}}
\def \R {{\bf {R}}}
\def \bv {v}

\def \bz {z}
\def \bx {x}
\def \by {y}

\def \H2s {H^{s+1}_0(\partial M\times [0,T/2])}

\def \diam {\hbox{diam }}

\def \e {\varepsilon}

\def \pa0 {\partial _0}

\def \p {\partial}

\def\e{\varepsilon}

\def\M{{M}}

\def \hat {\widehat }


\def\cM{{\M}}

\def\tilde{\widetilde}

\def\D{{\cal D}}

\def \mbeq {\begin {eqnarray}}
\def \meeq {\end {eqnarray}}

\def\la{\lambda}


\def \bfo {\begin {displaymath} }
\def \efo {\end {displaymath} }

\def \beq {\begin {eqnarray}}
\def \eeq {\end {eqnarray}}
\def \ba {\begin {eqnarray*}}
\def \ea {\end  {eqnarray*}}

\def \R {{\bf {R}}}
\def \bv { v}
\def \bu {{ {u}}}

\def \bz {{ {z}}}
\def \bx {{ {x}}}
\def \by {{ {y}}}

\def \H2s {H^{s+1}_0(\partial \M\times [0,T/2])}

\def \diam {\hbox{diam }}

\def \e {\varepsilon}

\def \pa0 {\partial _0}

\def \p {\partial}

\def \bK {\bf K}


\def \vol {\hbox{vol}}

\def \BSD {\hbox {boundary spectral data}}
\def \GH {\hbox {Gromov-Hausdorff}}
\def \bB {{\bf B}}
\def \bM {{\bf M}}
\def \diam {\hbox {diam}}
\def \inj {\hbox {inj}}


\def \tilde{\widetilde}

\title{Stability and Reconstruction
in Gel'fand Inverse Boundary Spectral Problem}

\author{Atsushi Katsuda\footnote{Department of Mathematics,
  Okayama University, Tsushima-naka,  Okayama, 700-8530, Japan }, 
Yaroslav Kurylev\footnote{
Department of Mathematical Sciences, Loughborough
University, Loughborough, LE11 3TU, UK}, and
  Matti Lassas\footnote{Rolf Nevanlinna Institute, University of Helsinki
Helsinki, P.O.Box 4, FIN-00014, Finland }
}
\maketitle

{\bf Abstract.} We consider stability and approximate reconstruction
of Riemannian manifold when the finite number of 
eigenvalues of the Laplace-Beltrami
operator and the boundary values of the corresponding eigenfunctions are given.
The reconstruction can be done in stable way when manifold is a priori
known to satisfy natural geometrical conditions related to curvature and
other invariant quantities.                    
\bigskip

{\bf I. Introduction} In this paper we consider the
questions of stability and approximate reconstruction
in the inverse boundary spectral problem
{\newtext which is also called the generalized Gelfand inverse problem 
\cite{Ge}  }
for
Riemannian manifolds. To formulate the problem and the main results, we need to
introduce some basic
notations. We will denote by $(\M,g)$ an
  unknown, $m$-dimensional,  compact connected
Riemannian manifold  with a (smooth) metric $g$ and non-empty boundary $\p \M$.
 Our goal is to find $(\M,g)$ from boundary data.
The boundary $\p \M$ is itself an $(m-1)$-dimensional compact
differentiable
manifold. We note that we do not assume the knowledge of the metric 
$i^*(g)$ generated on $\p\M$ by the metric $g$ ,
where $i: \p \M \to \M$ is an embedding or even the corresponding area
element
$d S_g$.

Because  in the inverse problems
the boundary $\p \M$ is usually known, we will consider
a class $\bM = \bM_{\p \M}$ of compact, connected Riemannian manifolds which
have the same, i.e. diffeomorphic, boundary, $\p \M$.

Let $-\Delta _g$ be the Laplace operator on $(\M,g)$
with Neumann boundary condition. Denote by
$\{\la_j, \, \varphi_j; \,j=1,2,\dots \}$
the complete set of eigenvalues and corresponding orthonormal
eigenfunctions
of $-\Delta _g,\,$ $0 =\la_1<\la_2,\dots,\, \varphi_1 = \vol ^{-1/2}
(\M,g)$,
where $\vol$ stands for the volume of $(\M,g)$.

Denote by $\bB $ the set 
of sequences $ \{\mu_j,\, \psi_j;\,$ $ j=1,2,\dots \}$ where
$\mu _j \in \R$ and $\psi_j \in L^2(\p \M)$
and by
$\D: \bM \to \bB $, the map
\bfo
\D\,(\M,g) = \left \{\la_j,\, \varphi _j|_{\p \M},\, j=1,2,\dots \right
\}.
\efo
\begin{definition}\label{def:1} 
{\newtext The collection 
  $ \left \{ \la_j,\, \varphi_j|_{\p \M}; \, j =1, 2, \dots \right \}$
is called {\it the boundary spectral data} of $(\M,g)$ and the map $\D$ -
 {\it the boundary spectral map}}.
\end{definition}

It was shown by Belishev-Kurylev \cite {BeKu} 
{\newtext who used the boundary control method \cite{Be} and the 
unique continuation result of Tataru \cite{Ta} }
that

\begin{theorem}\label{th:1} The map $\D: \bM \to \bB$ is injective.
\end{theorem}

{\newtext (For earlier uniqueness results in the Gelfand inverse
problem for isotropic operators obtained by the boundary control method see,
e.g. \cite{BeKu1}). The complex geometric optics method \cite{SU}
was applied to this problem in \cite{NSU}, \cite{No}.) }

In this paper, we will first analyse the question of stability of
the inverse problem,
i.e. the question of continuity of $\D ^{-1}$. Later,
we will also describe a procedure
of an approximate reconstruction of $(\M,g)$.

{\newtext Alessandrini \cite{Al},\cite{Al2}
 considered an inverse boundary problem for 
a Schr\"od\-ing\-er operator, $-\Delta + q$ in a bounded domain of $\R^m$
 and  obtained 
 a log-type stability estimate in the Gel'fand inverse 
problem (see also \cite{KuSt}).

Alessandrini and Sylvester \cite{AlS} and Stefanov and Uhlmann 
\cite{StU} considered
inverse boundary value problems for the wave equation which are  
closely related
to the Gel'fand inverse problem. Let
\begin{eqnarray*}
u_{tt} + a(\bx,D)u = 0, \quad \hbox{in} \quad  \Omega \times [0,T],
\quad u|_{t<0} < 0, \quad u|_{\p \Omega \times [0,T]} = f,
\end{eqnarray*}
where $a(\bx,D) = -\Delta + q$ in \cite{AlS} and 
$a(\bx,D) = - \p_ia^{ij}\p_j$ in \cite{StU} and the inverse data is given
as a non-stationary Dirichlet-Neumann map. For sufficiently large $T$ and,
in \cite{StU} $a^{ij}$ close to $\delta^{ij}$, they showed H\"older-type 
stability. These results raise the question of an optimal stability estimates
in the Gel'fand inverse problem and inverse boundary value problem for the
wave equation. Although exact stability estimates lie outside the scope
of this paper which is devoted to the analysis of geometric conditions
upon $(\M,g)$ to  provide continuity of $\D^{-1}$, we believe that 
H\"older-type estimates could be valid only under rather strong requirements 
upon $(\M,g)$, like  strong geodesic property (see e.g. \cite{Cr}, \cite{Sh}).
Further results on stability for different types of inverse boundary value 
problems could be found in \cite{Is}. For counterexamples see e.g.
\cite{M}. }

Clearly, to discuss the continuity of $\D ^{-1}$ it is necessary to
define appropriate topologies in $\bM$ and $\bB $. However,
it is well-known (see e.g. \cite {CoKr}) that inverse problems are,
in general, ill-posed. Therefore, to gain stability we need to
know {\it a priori} that the unknown object, which we are going
to reconstruct,  belongs
to some compact class $\bK$ (see e.g. \cite {Is}).
To obtain stability results, the topologies of $\bM$ and $\bB$
have to satisfy two properties. First, the map
$$\D: \bK \to \bB$$
has to be continuous. Second, $\D$ has to be injective. When
these two properties are satisfied and $\bK\subset \bM$
is a compact set
then the continuity
of the map
$$\D^{-1}: \D (\bK) \to \bK$$
follows
from basic results of  general topology.
In the first part of
the paper we carry out the constructions of the appropriate
topologies. We would
like to note that a  proper topology on $\bM$ is the Gromov-Hausdorff
topology and the (pre)-compact subsets we actually use are
the manifolds of bounded geometry
(see e.g.  \cite{BuBuI}, \cite {Gr}, \cite {Pe}, \cite{Sa}).
The topology on $\bB$ is defined from the 
boundary spectral convergence.

In the second part of the paper our attention will be focused
on the problem of an approximate reconstruction of $(\M,g)$ when we do not
know all $\BSD$. Rather, we know only the first eigenvalues,
$\la_j < \delta^{-1}$
{\newtext  with some small $\delta > 0$ }
and boundary values
of the corresponding eigenfunctions and, furthermore, we know them with
some
error. Namely, we have a finite collection of pairs
\bfo
\{\mu_j,\psi_j;\,j=1,\dots,{\it n}(\delta^{-1})\}, \quad \mu_j \in \R,\,
\psi_j \in L^2(\p \M)
\efo
and know that they are close to
$ \{ \la_j,\, \varphi_j|_{\p \M};\, j=1,\dots,{\it n}(\delta^{-1})\}$,
where ${\it n}(\delta^{-1})$ is the number of eigenvalues that 
are smaller than $\delta^{-1}$.
Our definition of topology on $\bB$ is adjusted to make the intuitive
notion
of closedness of $\{\mu_j,\psi_j\}$ to $\{ \la_j,\, \varphi_j|_{\p \M};\,
j\leq {\it n}(\delta^{-1})\}$, rigorous.

In the second part of the paper we describe the stable reconstruction
procedure but omit the proof showing that this procedure
is stable. These proof are to be published elsewhere.
We note that to carry out the reconstruction procedure
we will need to work under stronger requirements on $ (\M,g)$
than those that are  needed for stability.
Similarly, in the first part we will only formulate necessary
results omitting all proofs. The proofs and detailed description of the
reconstruction procedure will be given in a forthcoming paper.

\medskip

{\bf II. Main Results. Stability.}
We start with a proper topology on $\bM$.

\begin{definition}\label{def:1b} (Gromov-Hausdorff topology).
  Let $\e >0.$ The
Riemannian manifolds
$(\M^i,g^i) \in \bM,\,i=1,2$ are $\e$-close in the Gromov-Hausdorff
topology,
\bfo
d_{GH}((\M^1,g^1),\,(\M^2,g^2)) \leq  \e,
\efo
if there are
$\e$-nets $\{\bx_j^i;\,j=1,\dots,J(\e) \} \subset \M^i,\,
i=1,2,\,$
such that
\bfo
|d_{g_1}(\bx_j^1,\bx_k^1)-d_{g_2}(\bx_j^2,\bx_k^2)|\leq
\e, \quad j,k=1,\dots,J(\e),
\efo
where $d_{g^i}(\cdot,\cdot)$ stands for the distance on $(\M^i,g^i)$.
\end{definition}

{\newtext We remind that $\{\bx_j;\,j=1,\dots,J(\e) \} \subset \M$ is an
$\e$-net if
for every $\bx\in \M$ there is $\bx_j$
such that $d_{g}(\bx_j,\bx)\leq \e$. }


\medskip
\noindent
{\bf Remark 1.} In the future we will identify diffeomorphic manifolds
assuming implicitly that desired statements are valid after an
automorphism of $\M$.

\begin{definition}\label{def:2} 
(Riemannian manifolds of bounded geometry). For any $\Lambda,$  $D$, $i_0
>0$,
${\bM}(\Lambda, D, i_0 ) \subset \bM$ consists of Riemannian manifolds $(\M,g) \in \bM$
such that

\begin{tabular}{ll}
i) $\|\nabla Rm(\M,g)\|+\|Rm(\M,g)\| \leq \Lambda,$  & iii) $\diam (\M,g) \leq D$,\\ 
ii) $\|\nabla S(\M,g)\|+\|S(\M,g)\| \leq \Lambda, $ & iv) $\inj (\M,g) \geq i_0$.\\
\end{tabular}
\end{definition}

Here $Rm$ is 
{\newtext the Riemannian curvature tensor on $(\M,g)$ considered as a 4-linear 
form on $T\M$}, {$\nabla$ is the covariant derivative, 
 and}
 $S$ is the second
fundamental
form with respect to the inner product,
i.e. a quadratic form on $T \p \M$. Furthermore,
 $\diam$ is the diameter of $(\M,g)$ and $\inj$ is the
injectivity radius of $(\M,g)$.
Condition iv)
consists of three different
conditions:

\smallskip
\noindent
$(\alpha)$ Let $\bx \in \p \M$. Denote, as usual, by $B_{\p \M}(\bx,r)$
the
open metric ball (on $\p \M$) of the radius $r$ with center in $\bx$.
Then, for any $r < i_0$, $B_{\p \M}(\bx,r)$ is a domain
of normal coordinates on $\p \M$ centered at  $\bx$.

\smallskip
\noindent
$(\beta)$ Let $\bx \in \M$, $d(\bx,\p\M)\geq i_0$. Then, for
any
$r < i_0\,$, $B_{ \M}(\bx,r)$ is a domain
of normal coordinates on $ \M$ centered at  $\bx$.

\smallskip
\noindent
($\gamma$) Let $\bx \in \M$, $d(\bx,\p\M)<i_0$. Then,
  for any
$r < i_0$, the cylinder $C(\bx,r)$ is a domain of boundary normal
coordinates.
Here $C(\bx,r)$ consists of all points $\by \in \M$ such that
$d(\by,\p \M) < r$
and the unique boundary point $\bz(\by) \in \p \M,\,$ nearest to $\by$
is in the ball
$ B_{\p \M}(\bz(\bx),r)$.

\medskip
Next we introduce topology on the set ${\bB}$. 

\begin{definition}\label{def:3}
 (Boundary spectral topology.)
Let $\delta >0$.
Collections $$ \{\mu_j^i,\, \psi_j^i ;\, j=1,2,\dots\}\in \bB, \quad
i=1,2$$ are $\delta$-close if there are disjoint intervals
$$I_p = (a_p,b_p) \subset (-\delta ^{-1} -\delta,\ \delta^{-1}+\delta), 
\quad p=1,\dots, P,$$
such that

\smallskip
\noindent
i) $b_p-a_p <\delta,$

\smallskip
\noindent
ii) For any $\mu_j^i,\, i=1,2$ with  $\mu_j^i < \delta^{-1}$ 
there is $p$ such that
$\mu_j^i \in  I_p$.

\smallskip
\noindent
iii) For any $p$, the total number $n_p^i\, $ of $\mu_j^i$ inside $  I_p$
coincide, i.e.
$$ n_p^1 = n_p^2 \,(=n_p),$$

\smallskip
\noindent
iv) For any $p$ there is a unitary matrix $U_p = [u_p^{kl}]_{k,l=1}^{n_p}$
such that
$$
\|U_p \Psi_p^1 - \Psi_p^2\|_{L^2(\p \M)^{n_p}} \leq \delta.
$$
Here, $\Psi_p^i$ is the vector-function
$(\psi_{j(1)}^i \cdots \psi_{j(n_p)}^i)$, $j(1) < \cdots <j(n_p),$ with
$\mu_{j(k)} \in I_p$ for $k=1,\dots,n_p$.
\end{definition}

Clearly, condition iv) depends on the choice of the boundary measure
used in the definition of $L^2(\p \M)$. However, due to the compactness
of $\p \M$ different smooth measures on $\p \M$ determine
equivalent $L^2$-norms.

 We note that such topology was introduced by Alessandrini 
\cite{Al}, \cite{Al2} who studied 
stability in the
Gel'fand inverse problem for a Schr\"odinger operator.

\medskip
We are now able to formulate the principal stability
  result
for inverse boundary spectral problem.

\begin{theorem} \label{th:2} For any $\Lambda, D, i_0>0,$
  the map $\D^{-1}$ exists and is continuous on
$\bM (\Lambda, D,i_0)$. That is, 
there is $\delta >0$ such that if $(\M^i,g^i) \in {\bM}(\Lambda, D,
i_0), \,
i=1,2,$ and their \BSD,
\ba
  \left \{ \la_j^1,\, \varphi_j^1|_{\p \M}; \,
j =1, 2, \dots\right \}\quad\hbox{and}\quad
\quad \left \{ \la_j^2,\, \varphi_j^2|_{\p \M}; \,
j =1, 2, \dots\right \},
\ea
are $\delta$-close then
 $\M^1$ and $\M^2$ are diffeomorphic.
\noindent
 Moreover, for any $\alpha \in (0,1)$, 
$g^1 \to g^2$ in $C^{1,\alpha}$ -topology when $\delta \to 0$.
\end{theorem}

\medskip
\noindent
{\bf Remark 2.} {\newtext  The convergence of the metrics in the
$C^{1,\alpha}$-topology 
means that 
in a proper coordinate system on $\M =\M_i,\, i=1,2,$
the metric tensors $g_{ij}$ converge in the $C^{1,\alpha}$-sense,
see e.g. \cite{Pe}.}
%

\medskip
\noindent
{\bf Remark 3.} It follows from definition \ref{def:3} that to apply Theorem
\ref{th:2}
it is sufficient to know only the parts of the $\BSD$ of $(\M^1,g^1)$
and $(\M^2,g^2)$ that correspond to the eigenvalues $\la_j^i < \delta^{-1}$.

\smallskip
Theorem \ref{th:2} and the fact that in the class
$ {\bM}(\Lambda, D, i_0)$ the Gromov-Hausdorff topology
is equivalent to the Lipschitz topology (see (\cite {Kod})
implies the following corollary that describes the Lipschitz
convergence of Riemannian manifolds.

\begin{corollary}\label{cor:3} Let $(\M^i,g^i) \in {\bM}(\Lambda, D, i_0),
\,
i=1,2$. Then, for any $\e >0$ there is $\delta >0$ such that if the
$\BSD$ of these manifolds are $\delta$-close then
$\M^1$ and $\M^2$ are diffeomorphic and
$$
(1+\e)^{-1} \leq \frac{d_{g^1}(\bx,\by)}{d_{g^2}(\bx,\by)} \leq (1+\e), \,
\quad \bx,\by \in
\M.
$$
\end{corollary}

\medskip
\noindent
{\bf Remark 4.} It should be noted that ${\bM}(\Lambda, D, i_0)$ is
not closed in $\bM$ in the Gromov-Hausdorff topology. Therefore,
to obtain the desired stability result we should first extend Theorem
\ref{th:1} onto the closure
$\overline {{\bM}(\Lambda, D, i_0)}$. The set
$ {{\bM}(\Lambda, D, i_0)}$
consists of {some class of $C^{1,\alpha}$-smooth} Riemannian manifolds.
Although the boundary control method that was
instrumental to prove Theorem \ref{th:1} is, in general, not applicable
to  $C^{1,\alpha}$ Riemannian manifolds it was possible to show that
$$ \D:  \overline {{{\bM}(\Lambda, D, i_0)}} \to {{\bB}}$$
is injective (see section IV).

\medskip
 In Theorem \ref{th:2}, the condition iii) in ${\bM}(\Lambda, D, i_0)$ is 
automatically satisfied from an estimate of $\lambda_2$ in \cite{Me}. 
The other conditions  are
natural for the stability of the
inverse problem because, otherwise, the convergence of the $\BSD$
does not imply the convergence of the corresponding Riemannian manifolds. We illustrate
this
thesis with the following  examples.

\smallskip
\noindent
{\bf Example 1.}  Let $(\M,g)$ be a two-dimensional manifold obtained by the 
following surgeries of the two-sphere $S^2$. Add a handle $H$ near 
the North pole 
and cut-off 
a small piece $B$ to make a boundary near the South pole. 
If the metric size 
of the handle 
tend to $0$ uniformly
then, in the boundary spectral topology,
the $\BSD$  of $(\M,g)$ 
tend to the $\BSD$ 
 of the sphere with a hole but  without a handle  (see e.g. \cite{Tak}). 
This shows the necessities of i) 
and iv) in  
${\bM}(\Lambda, D, i_0)$. 


%

\begin{figure}[htbp]
\begin{center}
\psfrag{1}{$\partial M_1$}
\psfrag{2}{$\partial M_2$}
\psfrag{3}{$\partial M_1$}
\psfrag{4}{$\partial M_2$}
\psfrag{5}{$A$}
\includegraphics[width=12cm]{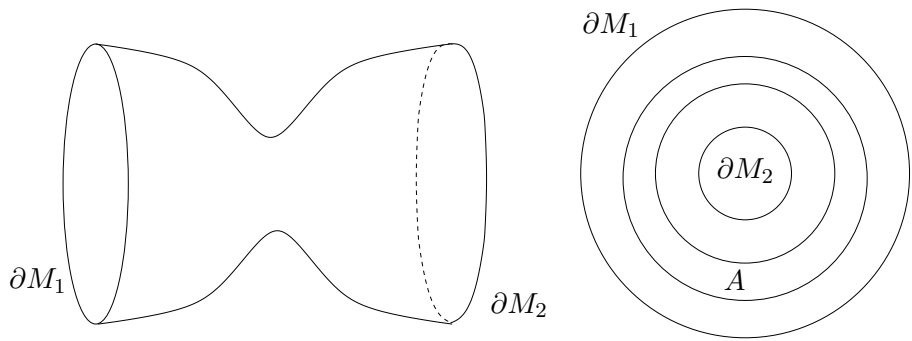}
\end{center}
\caption{Two representations of the same manifold:
Cylinder in ${\bf R}^3$ pinched in the middle
 and its representation as an annulus  in 
${\bf R}^2$
with variable metric. When $\e \to 0$ the arclength of circles in the region
$A$ tend to $0$ so that the manifold 
splits into two components. }
\end{figure}

\smallskip
\noindent
{\bf Example 2.} Let $(\M,g)$ be cylinder in ${\bf R}^3$
with the boundary made of two components,
$\p \M = \p \M_1 \cap \p \M_2$ (see Fig.1).
 We pinch the cylinder in the middle
so that the diameter tends
to $0$. 
 Therefore, the manifold eventually ``splits'' into two
disjoint semispheres. Meanwhile, 
in the boundary spectral topology on $\p \M$,
the $\BSD$ of $(\M,g)$ tend to the $\BSD$ of this disconnected manifold.
Again,
conclusion of the theorem is no more valid. 

 This second example may be also described in terms of a
 perturbation of
the metric in an annulus $A(1,4)= \{(r,\theta) : 1 \leq  r  \leq 4 \}\ ,$
where $(r,\theta)$) are the  polar coordinates in ${\bf R}^2$.
We define  the length element by
\bfo
dl^2_{\e} = dr^2 + \frac{r^2}{1 +\e^{-1} \chi(r)} \,d\theta^2, \quad
\chi(r) \in C^{\infty}_0,
\efo
where
\bfo
 \chi = 0 \quad \hbox{for} \quad r<2 \quad
\hbox{or} \quad
r>3, \quad
\chi >0 \quad \hbox{for} \quad r \in (2,3).
\efo
Then, for $\e \to 0$ each circle $r \in (2,3)$ shrinks to a point so that
the internal part near $r=1$ of $A$ becomes separated from the 
external part
near $r=4$.

We would note that although this example does not precisely satisfies
the definition of ${\bM}_{\p \M}$ because the limiting manifold is no longer
connected. However, this definition can be extended onto the case of, 
probably disconnected
Riemannian manifolds
with the same boundary $\p \M$ which have no components without the boundary
if we change the condition of boundedness of $\diam (\M)$ to the boundedness 
of $\vol (\M)$.

\smallskip
\noindent
{\bf Example 3} Take $X$ and $Y$ to be isospectral manifolds 
without boundary of
different 
topological types (see e.g. \cite {Sunada}) and $Z$ be a Riemannian
manifold with 
boundary $\p Z$.  Let $M_{\e}$ be a  Riemannian manifold obtained by
connecting $Z$ 
with $X$ by a thin tube of radius $\e ^{2}$ of  length $\e$. 
Similarly, let $N_{\e}$ be a  Riemannian manifold obtained by
connecting $Z$ 
with $Y$ by a thin tube of radius $\e ^{2}$ of  length $\e$. 
When $\e \to 0$, 
the boundary spectral data of $M_{\e}$ tend to 
the union of the boundary
spectral data of $Z$ and 
$\{\mu_j, 0\}$, where $\mu_j$ is the spectrum of $X$.
Clearly, the boundary spectral data of $N_{\e}$ tend to the same limit.
Thus, for
any $\delta > 0$, there exists $\e _0 > 0$ such that the 
boundary spectral data of 
$M_{\e}$ and $N_{\e}$ are $\delta$-close 
when $\e \leq \e _0$. However, $M_{\e}$ and $N_{\e}$ are 
of different topological types.



In these examples  $\|Rm\| \to \infty$ and 
$\hbox{inj} \to 0$  when $\e \to 0$ which
show the importance of conditions i) and iv) of the
theorem.  Similar examples can be given to 
illustrate the importance of condition ii).

\medskip
{\bf III. Main results. Approximate reconstruction.}
Definition \ref{def:1b} of the Gromov-Hausdorff topology in $\bM$ can be
easily extended to the
set of all compact metric spaces $(X,d)$. Thus, two Riemannian manifolds
$(\M^i,g^i)$ are $\e$-close in
the Gromov-Hausdorff metric if they possess $\e$-nets $X^i$ which,
being considered
as finite metric spaces with distance inherited from $(\M^i,g^i)$,
are $\e$-close
in the Gromov-Hausdorff topology. Therefore, to approximately reconstruct
the Riemannian manifold  $(\M,g)$ we should construct a
   metric space $(Y,d), Y = \{\by_1,\dots,\by_J\}$
which is $\e$-close to  $(\M,g)$, i.e.
\bfo
d_{GH}((\M,g),\,(Y,d)) \leq \e.
\efo

To construct a metric space approximation to $(\M,g)$ we will assume
further regularity properties of the class of manifolds.

%
%
%
%

\begin{theorem}\label{th:4.1} Let 
 $\Lambda, D,i_0 >0$. Then for
any $\e >0$ there is $\delta >0$ such that if
$ \{\mu_j,\, \psi_j|_{\p M};\,  j=1,2,\dots \} \in \bB \,$ is
$\delta$-close to the $\BSD$
$ \{ \la_j,\, \varphi_j|_{\p \M}; \,  j=1,2,\dots \}$
of a Riemannian manifold,
 $ (\M,g) \in \overline{{\bM}(\Lambda, D, i_0)} \, $,
then there is a finite metric space $(Y,d)$ such that
\bfo
d_{GH}((\M,g),\,(Y,d)) \leq \e
\efo
Moreover, there is a constructive algorithm to find $(Y,d)$
from
$ \{\mu_j^i,\, \psi_j^i;\,j=1,2,\dots \} \in \bB \,$.
\end{theorem}
The algorithm to construct $(Y,d)$ is described in  section V.

\medskip
\noindent
{\bf Remark 5.} By definition \ref{def:3} the $\delta$-closedness does not
involve
$ \{\mu_j^i,\, \psi_j^i  \}$ with $|\mu _j| >\delta^{-1}$. Similarly,
the construction of $(Y,d)$ does not use these data.

\medskip
{\bf IV. Geometric convergence and stability.}
The proof of Theorem \ref{th:2} is based upon the following result.

\begin{proposition}\label{th:5} For any $\Lambda, D,i_0 >0,\,$ the set
${\bM}(\Lambda, D, i_0)$ is pre-compact in the $\GH $ topology. Its
closure,
$\overline{{\bM}(\Lambda, D, i_0)}$ consists of differentiable manifolds
with $C^{1,\alpha}$-smooth metric, where $\alpha \in (0,1)$ is arbitrary.
\end{proposition}

If a sequence of Riemannian manifolds  $(\M^n,g^n)$ converges to $(\M,g)$ in the $\GH$
topology on
$\overline{{\bM}(\Lambda, D, i_0)}$, there is $n_0 \geq 1$ such that

\smallskip
\noindent

\smallskip
\noindent
 { 
i) For $n \geq n_0$, there is a diffeormorphism
 $F_n:\M\to \M^n$, i.e. 
 $\M^n$ and $\M$ are diffeomorphic.

\smallskip
\noindent
ii) for any $\alpha \in (0,1)$, $F_n^*(g^n)$ converge to  $g$ in 
the $C^{1,\alpha}$-topology. 
}
\smallskip
\noindent

\noindent
{\bf Remark 6.} { Proposition \ref{th:5} is obtained
by applying  the ideas developed by M.
Anderson (\cite {And}, see also \cite {HH}) in the interior of
the manifold and 
the fundamental equation of Riemannian geometry (see \cite{Pe})
near the boundary.
The latter equations is also known as the    
Riccati equation for the second fundamental form along normal 
geodesics. This equation gives the desired regularity estimates
for the metric tensor in the boundary normal coordinates.}
\smallskip
\noindent

The geometric convergence described by Proposition \ref{th:5} yields the
continuity
of the direct problem.

\smallskip
\begin{proposition}\label{p:6}
$
\D: \overline{{\bM}(\Lambda, D, i_0)} \to  {{\bB}}
$
is continuous.
\end{proposition}

\smallskip
\noindent
{\bf Remark 7.}
{Proposition \ref{th:5}
  was proven by Kodani \cite {Kod} in the  case of weaker regularity.
Namely, Kodani showed that the manifolds satisfying
\noindent

\begin{tabular}{ll}
i) $\|Rm(\M,g)\| \leq \Lambda,$ \hspace*{3cm} & iii) $\diam (\M,g) \leq D$,\\ 
ii) $\|S(\M,g)\| \leq \Lambda, $ & iv) $\inj (\M,g) \geq i_0$.\\
\end{tabular}

\noindent
are pre-compact in the Lipschitz topology rather
than the
$C^{1,\alpha}$-topology.}
This
result is sufficient for continuity of $\D$.
This can be obtained by means of  the perturbation 
theory of quadratic form 
(see e.g. \cite{Fu} in the case of
 manifolds without boundary). 
Therefore, the result of Kodani would have been sufficient to
prove
the continuity of $\D ^{-1} : \D(\overline {{\bM}(\Lambda, D, i_0)}) \to
\overline {{\bM}(\Lambda, D, i_0)}$
   if it were possible to show the
injectivity of $\D$ on $\overline {{\bM}(\Lambda, D, i_0)}$.
Unfortunately,
the method used to prove Theorem \ref{th:1} fails in
this case in its parts related to the approximate controllability
result by Tataru \cite {Ta} and also to the injectivity of the boundary
distance map (see the next section for the definition of this map).

Proposition \ref{th:5}, on the contrary, guarantees a stronger regularity of
the Riemannian manifolds
in $\overline {{\bM}(\Lambda, D, i_0)}$. It was used to extend
Theorem \ref{th:1} on the class
\bfo
 \hat{{\bM}} = \bigcup _{\Lambda,D,i_0} \overline {{\bM}(\Lambda, D,
i_0)}.
\efo

\begin{theorem}\label{th:1'} The map
$
\D : \hat{{\bM}} \to {{\bB}}$
is injective.
\end{theorem}

\medskip
{\bf V. Construction of a finite $\e$-net.}
The procedure of constructing an approximation $(Y,d)$
consists of several steps which we will explain separately. We will
also point out smoothness requirements which are  necessary
to carry out different
steps of the procedure.

\smallskip
{\bf a. Construction of Fourier coefficients of waves.}
For simplicity, we assume here that the metric on the boundary,
$i^*(g),\, i:\p\M \to \M$ is given.
  Consider an initial-boundary value
problem associated with the Neumann Laplace operator,
\bfo
u^f_{tt} - \Delta_gu^f = 0, \quad \hbox {in} \quad \M \times [0,D],
\efo
\bfo
\p_{\nu}u^f|_{\p \M \times [0,D]} = f \in C_p^1([0,D], \,L^2(\p \M)),
\efo
\beq
u^f|_{t=0} = u^f_t|_{t=0} =0,
\label {7}
\eeq
where $C_p^1([0,D], \,L^2(\p \M))$ consists of piecewise
$C^1$-functions of $t$. We
{\newtext  introduce the wave operator $W^t$ given by}
\bfo
W^t(f) = u^f(\cdot,t) \in C^1([0,D], \,L^2( \M))
\efo
and, if we denote by $u^f_k(t)$ the $k$-th Fourier coefficient of $u^f$,
\beq
u^f_k(t) = (u^f(t),\,\varphi_k) =
\int_0^t \,\int_{\p \M} s_k^t(\by,t')f(\by,t')\,dt'dS_g(\by),
\label {8}
\eeq
where $dS_g$ is the boundary area element and
\bfo
s_k^{t}(\by,t') = \frac{\sin (\sqrt {\la_k}(t-t'))}{\sqrt {\la_k}}\,
\varphi_k(\by).
\efo
This formula makes possible to find approximately the first Fourier
coefficients of the waves $u^f$ if we know
$ \{\mu_j, \psi_j;\, j=1,\dots,{\it n}(\delta^{-1})\}$.

\smallskip
{\bf b. Construction of domains of influence.}
Let $\Gamma \subset \p \M$ be open. Consider the waves $u^f(\tau)$ with
$f \in C_0^{\infty}(\Gamma \times [0,\tau]),\, \tau >0$. Clearly,
\bfo
\hbox {supp} [u^f(\tau)] \subset \M(\Gamma,\tau) =
\{\bx \in \M: d_g(\bx, \Gamma) \leq \tau \}.
\efo
By the fundamental result of Tataru (\cite {Ta},
see also \cite {KaKuLa}, Ch.2.5), we have

\smallskip
\begin{theorem}\label{th:7} Assume that the metric tensor $g$
is $C^1$-smooth.
For any $\tau >0$ and $\Gamma \subset \p \M$,
\bfo
\hbox {cl}_{L^2(\M)} \left \{u^f(\tau):
f \in C_0^{\infty}(\Gamma \times [0,\tau]) \right \} =
L^2(\M(\Gamma,\tau)),
\efo
where $L^2(\Omega) \subset L^2(\M)$ consists of functions with support in
$\overline {\Omega}$.
\end{theorem}
Let
$\eta >0$,  $D>\diam(\M)$ as in part iii) of Definition \ref{def:3}
 and
$
\Gamma_l,\quad l=1,\dots, L,$
be open subsets of $\p \M$ satisfying
\bfo
 \diam (\Gamma_l) < \eta,
\quad \Gamma_l \cap \Gamma_k = \emptyset, \quad
\p \cM = \overline{\cup \Gamma_l}.
\efo
 Let
\beq
\alpha = (\alpha_1,\cdots,
\alpha_L),\, \quad \alpha_l \leq D/\eta,
\label {8a}
\eeq
   be a multi-index. Denote by $\Sigma_{\alpha} \subset \p \cM \times [0,D]$
the
set
\bfo
\Sigma_{\alpha} = \bigcup _l \, (\Gamma_l \times [D-\alpha_l \eta,\,D])
\cap
(\p \cM  \times [0,D])
\efo
and by $\M_{\alpha}\subset \M$ the subdomain
\beq
\M_{\alpha} = \bigcup_l \M(\Gamma_l,\,\alpha_l\eta).
\label {3'}
\eeq

By Theorem \ref{th:7}, for any $\sigma >0$ there is $f=f_{\alpha}
\in C_0^{\infty}(\Sigma _{\alpha})$ such that
\bfo
\|u^{f}(D) - \chi_{\M_{\alpha}} \varphi_1 \|_{L^2(\M)}
\leq \sigma/2,
\efo
where $\chi_{\Omega}$ is the characteristic function of a set $\Omega$.

It then follows from the general results of control theory that there is a
finite linear combination $\hat {f}$,
\beq
\hat {f}(\bx,t) = \sum_{j=1}^J a_j s_j^D(\bx,t) \chi_{\Sigma_{\alpha}}(\bx,t),
\quad \bx \in \p\M,
\label {9a}
\eeq
such that
\beq
\|u^{\hat {f}}(D) - \chi_{\M_{\alpha}} \varphi_1 \|_{L^2(\M)}
\leq \sigma.
\label {10}
\eeq
The coefficients $a_j$ of the function (\ref{9a}) which satisfies
equation
(\ref{10}) can be found using the $\BSD$
$\{ \la_k,\, \varphi_k|_{\p \M};\,k=1,\dots,K\},\,$ $ K>J$ by means
of some variational
procedure.

To define this procedure rigorously, let $H_N(\Sigma_{\alpha}) \subset
C_p^1([0,D], \,L^2(\p \M))$ be an $N$-dimensional subspace spanned by the
functions $ \{s_n^D(\by,t)\,\chi_{\Sigma_{\alpha}} (\by,t);\,$
$n=1,\dots,N \}$.

Let $A_{N,I,K}(f)$ be a functional on $H_N(\Sigma_{\alpha})$,
\beq
A_{N,I,K}(f) =
\label {803}
\eeq
\bfo
= \sup \left \{ \left \|(P_KW^Df - \chi_{\M_{\alpha}}
\varphi_1,W^Dh) \right \|: \,
h \in H_I(\Sigma_{\alpha}), \quad \|P_KW^Dh \| \leq 1 \right \}.
\efo
Here $P_N$ is the orthoprojector in $L^2(\M)$,
\bfo
P_Nu = \sum_{n\leq N} (u,\varphi_n)\,\varphi_n.
\efo
Clearly, due to (\ref {8}) the functional $A(f)$ can be evaluated in terms
of
$ \{ \la_j,\, \varphi_j|_{\p \M};\,$ $j=1,\dots,{\it n}(\delta^{-1})\}$ if
${\it n}(\delta^{-1}) \geq \max (N,I,K)$. As we have $\{\mu_j,\, \psi_j
;\,j=1,\dots,
{\it n}(\delta^{-1})\}$ rather then
$ \{ \la_j,\, \varphi_j|_{\p \M};\,j=1,\dots,{\it n}(\delta^{-1})\}$, the
functional $A(f)$ can be evaluated only approximately.

\smallskip
\begin{theorem}\label{th:9} For any $\sigma >0$ there are parameters $C,
N,I,K,\delta,\,$ $N \leq I \leq K \leq {\it n}(\delta^{-1}),\,$ such that if
$ \left\{\mu_j,\, \psi_j
\right \}$ is $\delta$-close to the boundary spectral data
$\left \{ \la_j,\, \varphi_j|_{\p \M} \right \}$ of a Riemannian manifold 
$(\M,g)$,
then a minimizer $f^* \in H_N$ of the functional $A_{N,I,K}$ such that
\bfo
\|f^*\|_{L^2(\p\M \times [0,D])} \leq C,
\efo
satisfies the equation
\bfo
\left \|W^Df^* - \chi_{\M_{\alpha}} \varphi _1  \right \| \leq \sigma.
\efo
\end{theorem}

\medskip
By Proposition \ref{p:6} the $\BSD$ depend continuously
on the $\GH$
distance restricted on $\overline{{\bM}(\Lambda, D, i_0)}$.
Consider
\bfo
u^f(\cdot,t) = u^f_{(\M,g)}(\cdot,t)
\efo
as a function of $(\M,g) \in \overline{{\bM}(\Lambda, D, i_0)}$.
{\newtext By
Proposition \ref{th:5},  the
perturbation theory \cite{Ka}
shows that the wave $u^f_{(\M,g)}(\cdot,t)$
depends continuously (in the $C^1([0,D],L^2(\M))$-topology)
on the $\GH$ distance. As the set
$\overline{{\bM}(\Lambda, D, i_0)}$ is compact, the parameters}
$C,N,I,K$
and $\delta$ used in the above construction are uniform on
$\overline{{\bM}(\Lambda, D, i_0)}$ and depend only on $\sigma$.

\medskip
{\bf c. Construction of the boundary distance map.}
As
\bfo
\|\chi_{\M_{\alpha}} \varphi _1 \|^2 = \frac {\vol (\M_{\alpha})}{\vol
(\M)}, \quad \vol(\M) = \varphi_1(z)^{-2},\ z\in \p M,
\efo
Proposition \ref{th:5} makes it possible to evaluate an approximate volume,
$\vol ^a (\M_{\alpha})$. We will use this observation
to construct a finite approximation
to the set of the boundary distance functions associated with $(\M,g)$.
Namely, for any $\bx \in \M$, let $r_{\bx} \in L^{\infty}(\p\M)$ be given
by
\bfo
r_{\bx}(\bz) = d_g(\bx,\bz), \quad \hbox {for any} \quad \bz \in \p\M.
\efo
The boundary distance functions determine the boundary distance map
$R: \,(\M,g) \to L^{\infty}(\p\M)$,
\bfo
R(\bx) = r_{\bx}.
\efo
It turns out  that the map $R$ or, more precisely, its image
$R(\M,g)$ is sufficient to reconstruct $(\M,g)$ (see \cite {Ku}
for the procedure of the reconstruction of $(M,g)$ from $R(M)$).
We will use an approximation $R^*$ to $R(\M,g)$ to construct the desired
metric space $(Y,d)$ in Theorem \ref{th:4.1}.
Note that 
the metric $d$ is different to the 
induced metric as a (metric) subspace of $L^{\infty}(\p\M)$ (cf. 
Step 1 in subsection {\bf 
d}.

To find $R^*$ consider subdomains $\M_{\beta}^* \subset \M,$
\beq
\M_{\beta}^* = \left \{\bx \in \M: \, d_g(\bx, \Gamma _l) \in
((\beta _l -2) \eta, (\beta _l +2) \eta),\quad l=1,\dots,L \right \}.
\label {21}
\eeq
For any $\beta$ the subdomain $\M_{\beta}^*$ can be obtained as a finite
number of unions, intersections and compliments of the subdomains
$\M_{\alpha}$ of form (\ref{3'}). For any $\Omega, \Omega' \subset \M$,
\bfo
\vol (\Omega ^c) = \vol (\M) - \vol (\Omega),
\efo
where $\Omega ^c = \M \setminus \Omega$, and
\bfo
\vol (\Omega \cap \Omega') = \vol (\Omega) +\vol (\Omega') -
\vol (\Omega \cup \Omega').
\efo
Moreover, for any $\alpha, \beta$,
\bfo
\vol (\M_{\alpha} \cup \M_{\beta}) = \vol (\M_{\gamma}),
\quad \hbox {where} \quad \gamma _l = \max (\alpha _l, \beta _l).
\efo
Therefore, it is 
possible to evaluate an approximate volume $\vol ^a (\M_{\alpha}^*)$,
if we know $\{\mu_j,\, \psi_j\}$. Clearly, $| \vol ^a (\M_{\beta}^*)
-\vol  (\M_{\beta}^*) | < \sigma$ when  $\{\mu_j,\, \psi_j \}$
is $\delta$-close to
$\{ \la_j,\, \varphi_j|_{\p \M}  \}$ with sufficiently
 small $\delta$.

When $\vol ^a (\M_{\beta}^*) \geq \sigma$ we associate with this
$\beta$ a function
$r_{\beta} \in L^{\infty}(\p\M)$,
\bfo
r_{\beta}(\bz) = \beta _l \eta \quad \hbox {for} \, \quad \bz \in \Gamma
_l.
\efo
Thus, if $\vol (\M_{\beta}^*) \geq 2\sigma$ so that
$\vol ^a (\M_{\beta}^*) \geq \sigma$, there is an approximate function
$r_{\beta}$ such that for any $\bx \in \M_{\beta}^*$,
\bfo
\|r_{\beta} - r_{\bx} \| \leq 2 \eta.
\efo
{\newtext In the construction of the approximate volume
of a set $\M_{\beta}^* \neq \emptyset$ such that
 $\vol (\M_{\beta}^*) < 2\sigma$ 
it may happen that 
$\vol ^a (\M_{\beta}^*) < \sigma$.
Thus, according to the procedure, 
$R^*$ do not contain $r_{\beta}$  corresponding to this multi-index $\beta$.
Consider next a point $\bx \in \M_{\beta}^*$ for such $\beta$.
By condition (\ref {8a}) there is an upper bound 
$N(\eta)$ for the number of all multi-indexes $\beta$.
Moreover, due to 
condition iv) of
definition \ref{def:2}, $\vol (B_{\M}(\bx,\eta)) \geq c \eta ^m$. Therefore,
if $\sigma \leq C' \eta ^m N(\eta)^{-1}$, there is a multiindex
$\gamma, \ \gamma \neq \beta$ such that

\smallskip
\noindent
i. $\vol (\M_{\gamma}^*) \geq 2 \sigma$, so that
$\vol ^a (\M_{\gamma}^*) \geq  \sigma$,

\smallskip
\noindent
ii. $d(\bx,\M_{\gamma}^*) \leq \eta$.
\smallskip

\noindent
Thus 
\bfo
\|r_{\gamma} - r_{\bx} \| \leq 4 \eta.
\efo
Summarizing the previous considerations we obtain the following result.
}

\begin{lemma}\label{l:10} For any $(\M,g) \in \overline {\bM}$ and any
$\eta >0$ there
is  $\delta >0$ such that if $\{\mu_j,\, \psi_j\}$ is $\delta$-close
to the $\BSD\,$ $\{ \la_j,\, \varphi_j|_{\p \M}  \}$ then
the set $R^*$ constructed from $\{\mu_j,\, \psi_j;\,j=1,\dots,{\it
n}(\delta^{-1})\}$
is $\eta$-close to $R(\M,g)$, i.e.
\bfo
d_H(R^*,R(\M,g)) \leq \eta.
\efo
Here $d_H$ is the Hausdorff distance between subsets in
  the metric space $L^\infty(\p \M)$.
Moreover, $\delta$ can be chosen uniformly on any
$\overline{{\bM}(\Lambda, D, i_0)}$.
\end{lemma}
\medskip
{\bf d. Construction of an $\e$-net}
The set $R^*$ has a metric structure inherited from $L^{\infty}(\p\M)$,
\bfo
d_{\infty}(r_1,r_2) = \|r_1 -r_2\|_{L^{\infty}(\p\M)}.
\efo
If $(\M,g)$ is strongly  geodesic (e.g. \cite{Cr}, \cite{Sh}),
i.e. all geodesic inside $\M$ are minimal, then
\beq
d_g(\bx,\by) = \|r_{\bx} -r_{\by}\|_{L^{\infty}(\p\M)},
\label {23}
\eeq
for any $\bx,\by \in \M$. Thus,  the metric space $(R^*,d_{\infty})$ provides
a $2\eta$-net for $(\M,g)$, 
\bfo
d_{GH}((R^*,d_{\infty}),\,(\M,g)) \leq 2\eta.
\efo
However, condition (\ref{23}) is not valid for general manifolds.
Therefore, $d_{\infty}$-metric is inappropriate to construct
an approximation to $(\M,g)$ in the $\GH \,$ topology. In this section
we will describe how to equip
the finite metric space $R^*$ with another metric $d^*$ so that
$(Y,d)=(R^*,d^*)$ becomes close to $(\M,g)$ in the $\GH$ metric.
%
%
%
%
%
{\newtext In our description we will assume that the whole $R(\M)$
rather then $R^*$ is known. We will show how to find an approximation to
$d_g(\bx,\by)$ from $R(\M)$. This procedure is valid also in the case
when we know only $R^*$. The corresponding result is given at the end
of this section.}

\smallskip
\noindent
{\bf 1.} Assume first that $\bx_1,\bx_2 \in \M ^{int}$ can be connected by a
unique shortest geodesic, $\gamma (\bx_1,\bx_2)$. Assume in addition that
this geodesic can be continued beyond $\bx_2$ to
a boundary point $\bz \in \p\M$
so that $\gamma (\bx_1,\bz)$ is a shortest
geodesic between $\bx_1$ and $\bz$.  Then,
\bfo
d_g(\bx_1,\bx_2) = d_g(\bx_1,\bz) - d_g(\bx_2,\bz) =
\|r_{\bx_1} -r_{\bx_2}\|,
\efo
can be found from $R(\M,g)$.

It remains to identify those $r_1,r_2 \in R(M)$ which 
correspond to $\bx_1,\bx_2$
with this property.
{\newtext As $i^*g$ on $\p\M$ is known,}
it is
possible
to extend $(\M,g)$ to a larger Riemannian manifold 
$(\tilde {\M}, \tilde g)$ which has
a continuous, piecewise $C^1$ metric $ \tilde g$.  We can then extend each
function
$r_{\bx} \in C(\p\M)$  to the function
$\tilde r_{\bx} \in C(\tilde {\M} \setminus \M)$, $\,\tilde r_{\bx}(\by)=
d_{\tilde{\M}}(x,y)$.
{\newtext Take a point $\bz \in \p\M$ and consider minimal geodesics
$\gamma (\bx_i,\bz),\, i=1,2$. Assume that 
 $\gamma (\bx_i,\bz) \setminus \{\bz\}  \subset \M ^{{int}}$,
 $\bz$ is not a cut point along any of $\gamma (\bx_i,\bz)$ and
$\gamma (\bx_i,\bz)$ are transversal to $\p \M$.
Then these geodesics can be extended as minimal geodesics beyond $\bz$ 
upto some points $\by_i \in \tilde {\M} \setminus \M$.
 Denote by
$\gamma_i(\bz,\by_i) \subset \tilde {\M} \setminus \M$ these extensions.
Then the geodesic $\gamma(\bx_1,\bx_2)$ connecting $\bx_1$ and $\bx_2$
extends  to a minimal geodesic
to
$\bz \in \p \M$ if and only if $\gamma_1(\bz,
\by_1) = \gamma_2(\bz,
\by_2)$.}

\begin{figure}[htbp]
\begin{center}
\psfrag{1}{$y_3$}
\psfrag{2}{$x_2$}
\psfrag{3}{$\hspace{-3mm}x_1$}
\psfrag{4}{$\hspace{5mm}y_2$}
\psfrag{5}{$\gamma_{31}$}
\psfrag{6}{$\gamma_{32}$}
\psfrag{7}{$\gamma_{12}$}
\psfrag{8}{$\partial M$}
\psfrag{9}{$\partial \tilde M$}
\includegraphics[width=8cm]{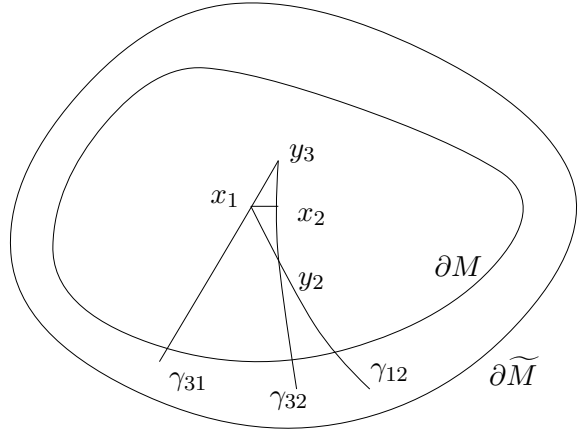} \label{pict3}
\end{center}
\caption{Example of a good triangle used to construct $d(x_1,x_2)$. In this 
case $x_1=y_1$ and $x_2$ lies on the geodesic $\gamma(y_2,y_3)$.}
\end{figure}

\smallskip
\noindent
{\bf 2.} Let now $\bx_1,\bx_2 \in \M ^{{int}}$ lie
on the sides $\gamma (\by_3,\by_1)$ and
$\gamma (\by_3,\by_2)$ of a
{\newtext good}
geodesic triangle $\Delta
(\by_1,\by_2,\by_3)$.
{\newtext A geodesic triangle $\Delta$ is called good if
}
all three geodesics $\gamma (\by_i,\by_j)$ can be continued
to $\p\M$ as the shortest geodesics (see Fig. 2).
   Then, by part {\bf 1} we can find $d_g(\by_i,\by_j) \,$ and also
$d_g(\bx_i,\by_i)$. Let
$\tilde{\Delta} (\tilde{\by}_1,\tilde{\by}_2,\tilde{\by}_3)$
be a triangle in $\R^2$ with $d_e(\tilde{\by}_i,\tilde{\by}_j) =
d_g(\by_i,\by_j)$, where $d_e$ stands for the Euclidean distance. Let
$\tilde{\bx}_i, i=1,2$ be the points on the sides
$[\tilde{\by}_3,,\tilde{\by}_i]$
of $\tilde{\Delta}$ with $d_e(\tilde{\bx}_i,\tilde{\by}_i) =
d_g(\bx_i,\by_i)$.
Then, by the Alexandrov lemma,
\bfo
|d_g(\bx_1,\bx_2) - d_e(\tilde{\bx}_1,\tilde{\bx}_2)| \leq c \sigma^2,
\efo
where $\sigma = \max d_g(\by_i,\by_j)$ and $c$ is uniform on
$\overline{{\bM}(\Lambda, D, i_0)}$.

Therefore, if for some given $\bx_1,\bx_2$,
we can find a proper geodesic triangle
$\Delta (\by_1,\by_2,\by_3)$ of a small size, we can find
$d_g(\bx_1,\bx_2)$
with an error of order $\sigma ^2$.

\smallskip
\noindent
{\bf 3.} The existence of a good geodesic triangle for any sufficiently
close
$\bx_1,\bx_2 \in \M ^{{int}}$ is based on the following result.

\begin{lemma}\label{l:11} Let $\Lambda,D,i_0 >0$ and
$(\M,g) \in
\overline{{\bM}(\Lambda, D, i_0)}$.
Let $\bx \in \M$.
 Then there is  a constant $\rho>0$ such that
for any
$(\M,g) \in
\overline{{\bM}(\Lambda, D, i_0)}$ and $\bx \in (\M,g)$
with $d_g(\bx,\p\M) \geq C \rho$ there is a vector
  $\bu \in T_{\bx}\M\,$ $|\bu|=1$ with the following property: 

\smallskip
\noindent
Let $(\by,\bv) \in B_{\rho}(\bx,\bu),\,|\bv| =1$. Then the geodesic
$\gamma _{\by}(t\bv),\, t\geq 0,\,$ can be extended as a shortest
geodesic to the boundary.

\end{lemma}
Here $B_{\rho}(\bx,\bu)$ is the ball of radius $\rho$ with center
at $(\bx,\bu)$ in the Sasakian metric on $T\M$.
{\newtext Combining steps {\bf 1-3} we can find an approximation
to $d_g(\bx_1,\bx_2)$ when $\bx_1,\bx_2$ are sufficiently
close to each other and not too close to $\p\M$.
(For the definition 
of the Sasakian 
metric, see e.g. \cite {Sa}.)

\smallskip
\noindent
{\bf 4.} Let
\bfo
d_g(\bx_i,\p\M) = \min_{\bz \in \p\M} r_i(\bz) \leq \eta < i_0,
\quad r_i =r_{x_i}.
\efo
There are unique $\bz_i \in \p\M$ with
$d_g(\bx_i,\bz_i) = r_i(\bz_i)$. Denote
\bfo
\tilde{d} = \tilde{d}(\bx_1,\bx_2) =
\left[d^2_{\p\M}(\bz_1,\bz_2) + |r_1(\bz_1) - r_2(\bz_2)|^2 \right]^{1/2},
\efo
where $d_{\p\M}(\cdot,\cdot)$ is the distance along $\p\M$. Then,
\bfo
|d_g(\bx_1,\bx_2) - \tilde{d}(\bx_1,\bx_2)| \leq C(\eta \tilde{d}
+\tilde{d}^2).
\efo
Therefore, $\tilde{d}$ is a good approximation to $d$ when $\bx_i$ are close to
$\p\M$.

\smallskip
\noindent
{\bf 5.} Steps {\bf 1-4} define approximate distance $ \tilde{d}(\bx_1,\bx_2)$
for sufficiently close $\bx_1,\,\bx_2$. Using a standard procedure
we can extend $ \tilde{d}(\bx_1,\bx_2)$ onto arbitrary $\bx_i \in \M$.
}

\smallskip
When we know the metric space $(R^*,d_{\infty})$ rather then $R(\M,g)$
we can carry out the above constructions approximately.
This gives rise to the following result.

\begin{lemma}\label{l:12}
Let $Y  \subset L^{\infty}(\p\M)$  be a set described
in Lemma \ref{l:10}.  It is possible to equip $Y$ with
a metric structure $d(r_1,r_2)$
 so that
\bfo
d_{GH}((Y,d),(\M,g)) \leq \e (\eta),
\efo
where $\e (\eta) \to 0\,$ when $\eta \to 0$. The function $\e(\eta)$
can be found uniformly for $\overline{{\bM}(\Lambda, D, i_0)}$.
\end{lemma}

\smallskip
Combining steps {\bf a.-d.} we obtain a procedure to construct a finite
$\e$-net for $(\M,g),\, \e = \e(\delta)$.

\medskip
{\bf VII. Concluding remarks.}

   {\bf Remark 8.} The above considerations lack quantitative
estimate for the function $\e(\delta)$. Combining our construction
with the Carleman estimates (see, e.g. \cite {Ta}, \cite {Ta1}, \cite {KaKuLa})
we will  obtain quantitative estimates for  $\e =\e(\delta)$.
This will be described in a forthcoming paper.

\smallskip
{\bf Remark 9.} All considerations remain valid for the Dirichlet
Laplacian. We need just to take into account that $\varphi_1 \neq 0$
in $\M^{{int}}$.

\smallskip
{\bf Remark 10.} It is possible to use the convergence of the
heat kernels restricted  to $\p\M$ instead of the 
boundary spectral 
convergence.
Besides, it is possible to modify the procedure to find an
approximation to the heat kernel associated with $-\Delta _g$ rather then
$R^*$. This puts our results
within the framework
of the spectral convergence of Riemannian manifolds (see e.g. 
\cite {BeBeG}, \cite {KaKu}). This will be described in the forthcoming paper.

\smallskip
{\bf Remark 11.} Requirements of $C^1$-regularity of $Rm$ and $S$ are of 
technical nature. We intend to improve reconstruction procedure to be valid 
for weaker regularity assumptions for  $Rm$ and $S$.

Part of our results also briefly announced in \cite {Kat}.

\medskip
{\bf Acknowledgements}
The authors would like to express their gratitude to Y. D. Burago,
D. Y. Burago,
M. Gromov, I. M. Gel'fand, A. Katchalov, T. Sakai, 
E. Somersalo, J. Takahashi
for helpful discussions and friendly support. This work was partly
supported
by EPSRC (UK) grant EPSRC (UK) grants GR/M14463 and 36595, RiP Program,
Oberwolfach (Germany), Grant-in Aid for Scientific Research 
09640109 and 12640073 (Japan), IHES (France),
 Finnish Academy project 42013 and 172434, MSRI NSF grant DMS-9810361 (USA),
 and TEKES (Finland).
  We are grateful to all these organizations.

\end {document}